\documentclass[a4paper,11pt]{article}

\usepackage[latin1]{inputenc}
\usepackage{amsmath}
\usepackage{amsthm}
\usepackage{amssymb}
\usepackage{mathrsfs}
\usepackage{graphicx}
\usepackage{hyperref}

\title{Hyperbolic diagram groups are free}
\date{\today}
\author{Anthony Genevois}

\begin{document}

\maketitle

\begin{abstract}
In this paper, we study the so-called \textit{diagram groups}. Our main result is that diagram groups are free if and only if they do not contain any subgroup isomorphic to $\mathbb{Z}^2$. As an immediate corollary, we get that hyperbolic diagram groups are necessarily free, answering a question of Guba and Sapir.
\end{abstract}

\tableofcontents

\section*{Introduction}
\addcontentsline{toc}{section}{Introduction}

\hspace{0.6cm} It is a well-known fact that a hyperbolic group does not contain any subgroup isomorphic to $\mathbb{Z}^2$. Although the converse does not hold in general, containing $\mathbb{Z}^2$ turns out to be the only obstruction to hyperbolicity for many classes of groups. For instance, Thuston's Hyperbolization Conjecture, proved by Perelman, implies that the fundamental group of a closed 3-manifold is hyperbolic if and only if it does not contain $\mathbb{Z}^2$. In the context of CAT(0) cube complexes, the following two results are known:

\medskip \noindent
\textbf{Theorem:} [CH09, Corollary 4] \textit{The fundamental group of a compact nonpositively curved special cube complex is hyperbolic if and only if it does not contain $\mathbb{Z}^2$.}

\medskip \noindent
\textbf{Theorem:} [SW11, Theorem 1.3] \textit{Let $X$ be a nonpositively curved cube complex whose universal cover satisfies the cyclic facing triples condition. Then its fundamental group $\pi_1(X)$ is hyperbolic if and only if it does not contain $\mathbb{Z}^2$.}

\medskip \noindent
In this paper, we study the so-called \textit{diagram groups}. The definition of diagram groups was first given by Meakin and Sapir, with the first results found by their student Kilibarda in her thesis [Kil94]. Although it was proved that diagram groups define a large class of groups with strong properties [GS97, GS99, GS06a, GS06b], very little is known on their geometric properties. Nevertheless, an important property due to Farley [Far00] is that they act freely on a CAT(0) cube complex. Therefore, a natural problem, to which the two theorems above do not apply, is to know whether or not a diagram group is hyperbolic if and only if it does not contain $\mathbb{Z}^2$.

In fact, Guba and Sapir [GS06a] gave a positive answer to this problem for a large class of diagram groups, namely for diagram groups \textit{over a complete directed 2-complex}. Using a cohomological argument, they actually proved much more: such a diagram group is free if and only if it does not contain any subgroup isomorphic to $\mathbb{Z}^2$.

In this article, we extend their result to any diagram group by exploiting the action of diagram groups on Farley's CAT(0) cube complexes. Explicitly, our main result is:

\medskip \noindent
\textbf{Main Theorem:} \textit{A diagram group is free if and only if it does not contain $\mathbb{Z}^2$. }

\medskip \noindent
It partially answers [GS06a, Problem 9.15] by extending [GS06a, theorem 9.14] to any diagram group. As an immediate consequence, hyperbolic diagram groups turn out to be free, answering a weaker question appearing in [S07, Problem 6.1]:

\medskip \noindent
\textbf{Corollary:} \textit{Let $G$ be a finitely generated diagram group. Then $G$ is hyperbolic if and only if $G$ does not contain any subgroup isomorphic to $\mathbb{Z}^2$ if and only if $G$ is free.}

\medskip \noindent
Although hyperbolic diagram groups turn out to be free, there exist non-free hyperbolic subgroups of diagram groups. More precisely, in [CSS07], an example of a diagram group (in fact a right-angled Artin group) with a hyperbolic surface subgroup is given. Characterizing hyperbolic subgroups of diagram groups remains a wide open problem.

\medskip \noindent
The plan of the paper is as follows. In Section 1, we give some preliminaries on diagram groups, CAT(0) cube complexes, and Squier's and Farley's complexes. In Section 2, we describe hyperplanes in Farley's CAT(0) cube complex and characterize their stabilizers in several cases. In Section 3, we prove several lemmas about hyperplanes in Squier's nonpositively curved cube complex. Finally, we prove our main theorem in Section 4, and we show in Section 5 how it can be used in practice to determine whether or not a given diagram group is free.

\medskip
Here is a sketch of the proof of our main theorem. Let $G$ be a diagram group without any subgroup isomorphic to $\mathbb{Z}^2$ and suppose that $G$ is not free. In particular, a non-trivial relation holds in $G$ with respect to some generating set $S$, and because such a relation uses only finitely many letters of $S$, it holds in some finitely generated subgroup $H$. Now, thanks to Grushko's theorem, $H$ can be written as the free product of finitely many freely irreducible groups; of course, at least one of these groups is not free, since otherwise $H$ would be free, so we may suppose without loss of generality that $H$ is freely irreducible. Therefore, to get a contradiction, it is sufficient to prove that any non-trivial finitely generated subgroup of $G$ either splits as a non-trivial free product or is infinite cyclic; equivalently, it is sufficient to prove that such a subgroup has to split non-trivially over the trivial subgroup.

Now let $G$ act on its Farley's CAT(0) cube complex $X$. We prove that it is possible to find a hyperplane $J \subset X$ satisfying the following two conditions:\\
\indent $\bullet$ $g_1 \cdot J$ and $g_2 \cdot J$ are either equal or disjoint, for all $g_1,g_2 \in G$,\\
\indent $\bullet$ the stabilizer $\mathrm{stab}(J)$ is trivial.\\
Thus, the orbit of hyperplanes $G \cdot J$ induces an \textit{arboreal structure} on $X$, i.e., if we define the graph $T$ whose vertices are the connected components of $X \backslash H \cdot J$ and where two components are linked by an edge precisely when they are separated by exactly one translate of $J$, then $T$ is a tree. Therefore, we get an action of $H$ on a simplicial tree $T$. Moreover, there is only one orbit of edges, whose stabilizer is trivial.

By a classical argument of Bass-Serre theory [Se80], in order to conclude that $H$ splits non-trivially over the trivial subgroup, it is sufficient to prove that the action $H \curvearrowright T$ has no global fixed point. We will finish the argument by proving that this is the case if $J$ is well-chosen.

\medskip \noindent
\textbf{Acknoledgements.} I would like to thank Michah Sageev for having suggested me to look at CAT(0) cube complexes associated to diagram groups, and of course my advisor, Peter Haïssinsky, for all our discussions.

\section{Preliminaries}

\noindent
\textbf{1.1. Diagram groups.} We refer to [GS97, §3 and §5] for a detailed introduction to \textit{semigroup diagrams} and \textit{diagram groups}. 

For an alphabet $\Sigma$, let $\Sigma^+$ denote the free semigroup over $\Sigma$. If $\mathcal{P}= \langle \Sigma \mid \mathcal{R} \rangle$ is a semigroup presentation, where $\mathcal{R}$ is a set of pairs of words in $\Sigma^+$, the semigroup associated to $\mathcal{P}$ is the one given by the factor-semigroup $\Sigma^+ / \sim$ where $\sim$ is the smallest equivalent relation on $\Sigma^+$ containing $\mathcal{R}$. For convenience, we will assume that if $u=v \in \mathcal{R}$ then $v=u \notin \mathcal{R}$; in particular, $u=u \notin \mathcal{R}$.

A \textit{semigroup diagram over $\mathcal{P}$} is the analogue for semigroups to van Kampen diagrams for group presentations. Formally, it is a finite connected planar graph $\Delta$ whose edges are oriented and labelled by the alphabet $\Sigma$, satisfying the following properties:\\
\indent $\bullet$ $\Delta$ has exactly one vertex-source $\iota$ (which has no incoming edges) and exactly one vertex-sink $\tau$ (which has no outgoing edges);\\
\indent $\bullet$ the boundary of each cell has the form $pq^{-1}$ where $p=q$ or $q=p \in \mathcal{R}$;\\
\indent $\bullet$ every vertex belongs to a positive path connecting $\iota$ and $\tau$;\\
\indent $\bullet$ every positive path in $\Delta$ is simple.\\
In particular, $\Delta$ is bounded by two positive paths: the \textit{top path}, denoted $\text{top}(\Delta)$, and the \textit{bottom path}, denoted $\text{bot}(\Delta)$. By extension, we also define $\text{top}(\Gamma)$ and $\text{bot}(\Gamma)$ for every \textit{subdiagram} $\Gamma$. In the following, the notations $\text{top}(\cdot)$ and $\text{bot}(\cdot)$ will refer to the paths and to their labels. Also, a \textit{$(u,v)$-cell} (resp. a \textit{$(u,v)$-diagram}) will refer to a cell (resp. a semigroup diagram) whose top path is labelled by $u$ and whose bottom path is labelled by $v$.

Two words $w_1,w_2$ in $\Sigma^+$ are \textit{equal modulo $\mathcal{P}$} if their images in the semigroup associated to $\mathcal{P}$ coincide. In particular, there exists a \textit{derivation} from $w_1$ to $w_2$, i.e., a sequence of relations of $\mathcal{R}$ allowing us to transform $w_1$ into $w_2$. To any such derivation is associated a semigroup diagram, or more precisely a $(w_1,w_2)$-diagram, whose construction is clear from the example below. As in the case for groups, the words $w_1,w_2$ are equal modulo $\mathcal{P}$ if and only if there exists a $(w_1,w_2)$-diagram.

\medskip \noindent
\textbf{Example 1.1.1:} Let $\mathcal{P}= \langle a,b,c \mid ab=ba, ac=ca, bc=cb \rangle$ be a presentation of the free Abelian semigroup of rank three. In particular, the words $a^2bc$ and $caba$ are equal modulo $\mathcal{P}$, with the following possible derivation:

\begin{center}
$aabc \overset{(a,ab \to ba,c)}{\longrightarrow} abac \overset{(ab,ac \to ca, \emptyset)}{\longrightarrow} abca \overset{(a,bc \to cb,a)}{\longrightarrow} acba \overset{(\emptyset,ac \to ca,ba)}{\longrightarrow} caba$.
\end{center}

\noindent
Then, the associated $(a^2bc,caba)$-diagram $\Delta$ is:

\begin{center}
\includegraphics[scale=0.6]{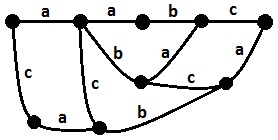}
\end{center}

\noindent
On such a graph, the edges are supposed oriented from left to right. Here, the diagram $\Delta$ has nine vertices, twelve edges and four cells; notice that the number of cells of a diagram corresponds to the length of the associated derivation. The paths $\text{top}(\Delta)$ and $\text{bot}(\Delta)$ are labelled respectively by $a^2bc$ and $caba$. 

Since we are only interested in the combinatorics of semigroup diagrams, we will not distinguish isotopic diagrams. For example, the two diagrams below will be considered as equal.

\begin{center}
\includegraphics[scale=0.6]{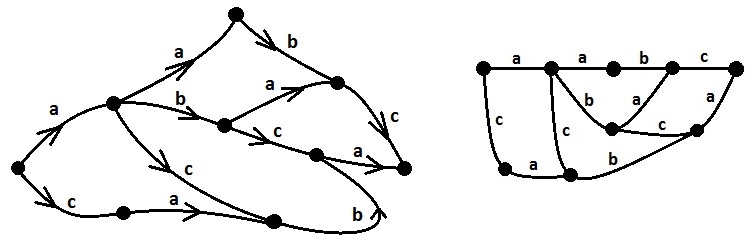}
\end{center}

\noindent
If $w \in \Sigma^+$, we define the \textit{trivial diagram} $\epsilon(w)$ as the semigroup diagram without cells whose top and bottom paths, labelled by $w$, coincide. Any diagram without cells is trivial. A diagram with exactly one cell is \textit{atomic}.

If $\Delta_1$ is a $(w_1,w_2)$-diagram and $\Delta_2$ a $(w_2,w_3)$-diagram, we define the \textit{concatenation} $\Delta_1 \circ \Delta_2$ as the semigroup diagram obtained by identifying the bottom path of $\Delta_1$ with the top path of $\Delta_2$. In particular, $\Delta_1 \circ \Delta_2$ is a $(w_1,w_3)$-diagram. Thus, $\circ$ defines a partial operation on the set of semigroup diagrams over $\mathcal{P}$. However, restricted to the subset of $(w,w)$-diagrams for some $w \in \Sigma^+$, it defines a semigroup operation; such diagrams are called \textit{spherical with base $w$}. We also define the \textit{sum} $\Delta_1+ \Delta_2$ of two diagrams $\Delta_1,\Delta_2$ as the diagram obtained by identifying the rightmost vertex of $\Delta_1$ with the leftmost vertex of $\Delta_2$.

\begin{center}
\includegraphics[scale=0.6]{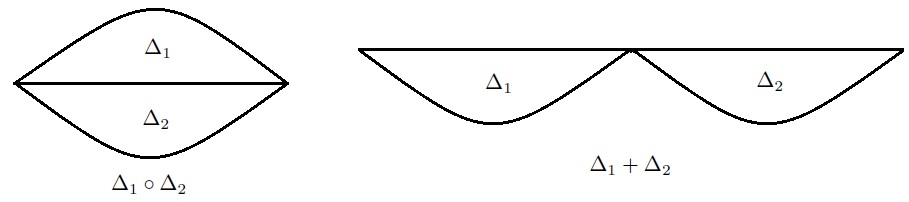}
\end{center}

\noindent
Notice that any semigroup diagram can be viewed as a concatenation of atomic diagrams. In the following, if $\Delta_1,\Delta_2$ are two diagrams, we will say that $\Delta_1$ is a \textit{prefix} (resp. a \textit{suffix}) of $\Delta_2$ if there exists a diagram $\Delta_3$ satisfying $\Delta_2= \Delta_1 \circ \Delta_3$ (resp. $\Delta_2= \Delta_3 \circ \Delta_1$). Throughout this paper, the fact that $\Delta$ is a prefix of $\Gamma$ will be denoted by $\Delta \leq \Gamma$.

Suppose that a diagram $\Delta$ contains a $(u,v)$-cell and a $(v,u)$-cell such that the top path of the first cell is the bottom path of the second cell. Then, we say that these two cells form a \textit{dipole}. In this case, we can remove these two cells by first removing their common path, and then identifying the bottom path of the first cell with the top path of the second cell; thus, we \textit{reduce the dipole}. A diagram is called \textit{reduced} if it does not contain dipoles. By reducing dipoles, a diagram can be transformed into a reduced diagram, and a result of Kilibarda [Kil94] proves that this reduced form is unique. If $\Delta_1,\Delta_2$ are two diagrams for which $\Delta_1 \circ \Delta_2$ is well defined, let us denote by $\Delta_1 \cdot \Delta_2$ the reduced form of $\Delta_1 \circ \Delta_2$.

If $w \in \Sigma^+$, we define the \textit{diagram group} $D(\mathcal{P},w)$ as the set of reduced $(w,w)$-diagrams endowed with the product $\cdot$ we defined above. If $\Delta$ is a $(w_1,w_2)$-diagram, let $\Delta^{-1}$ denote the $(w_2,w_1)$-diagram obtained from $\Delta$ by a mirror reflection with respect to $\mathrm{top}(\Delta)$. It can be noticed that, if $\Delta$ is a spherical diagram, then $\Delta^{-1}$ is the inverse of $\Delta$ in the associated diagram group.

Although this definition of $D(\mathcal{P},w)$ does not seem to give much information on its group structure, it allows to define a class of \textit{canonical subgroups}. If $\Gamma$ is a $(w,u)$-diagram and if we write $u=x_1u_1 \cdots x_nu_nx_{n+1}$, where the $x_i$ and $u_i$ are (possibly empty) subwords of $u$, then the map

\begin{center}
$(U_1,\ldots, U_n) \mapsto \Gamma \cdot \left( \epsilon(x_1) + U_1 + \cdots + \epsilon(x_n)+ U_n+ \epsilon(x_{n+1}) \right) \cdot \Gamma^{-1}$
\end{center}

\noindent
defines an embedding from $D(\mathcal{P},u_1) \times \cdots \times D(\mathcal{P},u_n)$ into $D(\mathcal{P},w)$.

\medskip \noindent
\textbf{1.2. CAT(0) cube complexes.} A \textit{cube complex} is a CW-complex constructed by gluing together cubes of arbitrary (finite) dimension by isometries along their faces. Furthermore, it is \textit{nonpositively curved} if the link of any of its vertices is \textit{flag} (ie. $n+1$ vertices span a $n$-simplex if and only if they are pairwise adjacent), and \textit{CAT(0)} if it is nonpositively curved and simply-connected. See [BH] for more information.

Although a CAT(0) cube complex $X$ can be endowed with a CAT(0) metric, it is often more convenient to introduce a more ``combinatorial" distance. We define the \textit{combinatorial distance} $d_c$ on the set of vertices $X^{(0)}$ of $X$ as the graph metric associated to the 1-skeleton $X^{(1)}$. In particular, \textit{combinatorial geodesics} are edge-paths of minimal length, and a subcomplex is \textit{combinatorially convex} if it contains any combinatorial geodesic between two of its points. In [Che00], Chepoi notices that the class of graphs appearing as 1-skeletons of CAT(0) cube complexes coincides with the class of \textit{median graphs}, which we now define.

Let $\Gamma$ be a graph. If $x,y,z \in \Gamma$ are three vertices, a vertex $m$ is called a \textit{median point of $x,y,z$} whenever

\begin{center}
$d(x,y)=d(x,m)+d(m,y)$, $d(x,z)=d(x,m)+d(m,z)$, $d(y,z)=d(y,m)+d(m,z)$.
\end{center}

\noindent
Notice that, for every geodesics $[x,m]$, $[y,m]$ and $[z,m]$, the concatenations $[x,m] \cup [m,y]$, $[x,m] \cup [m,z]$ and $[y,m] \cup [m,z]$ are also geodesics; furthermore, if $[x,y]$, $[y,z]$ and $[x,z]$ are geodesics, then any vertex of $[x,y] \cap [y,z] \cap [x,z]$ is a median point of $x,y,z$.

The graph $\Gamma$ is \textit{median} if every triple of pairwise distinct vertices admits a unique median point.

\medskip \noindent
\textbf{Theorem 1.2.1:} [Che00, Theorem 6.1] \textit{A graph is median if and only if it is the 1-skeleton of a CAT(0) cube complex.}

\medskip \noindent
Another important feature of cube complexes is the notion of \textit{hyperplane}. Let $X$ be a nonpositively curved cube complex. Formally, a \textit{hyperplane} $J$ is an equivalence class of edges, where two edges $e$ and $f$ are equivalent whenever there exists a sequence of edges $e=e_0,e_1,\ldots, e_{n-1},e_n=f$ where $e_i$ and $e_{i+1}$ are parallel sides of some square in $X$. Notice that a hyperplane is uniquely determined by one of its edges, so if $e \in J$ we say that $J$ is the \textit{hyperplane dual to $e$}. Geometrically, a hyperplane $J$ is rather thought of as the union of the \textit{midcubes} transverse to the edges belonging to $J$.

\begin{center}
\includegraphics[scale=0.4]{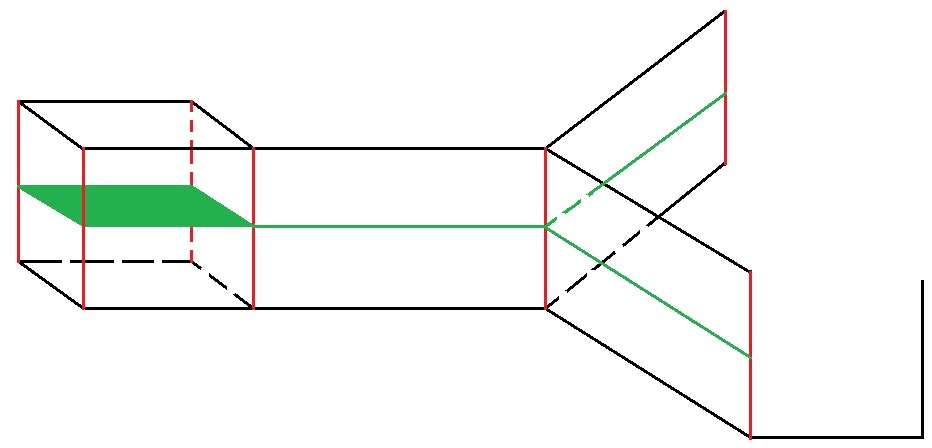}
\end{center}

\noindent
The \textit{neighborhood} $N(J)$ of a hyperplane $J$ is the smallest subcomplex of $X$ containing $J$, i.e., the union of the cubes intersecting $J$. In the following, $\partial N(J)$ will denote the union of the cubes of $X$ intersecting $N(J)$ but not $J$, and $X \backslash \backslash J= \left( X \backslash N(J) \right) \cup \partial N(J)$. Notice that $N(J)$ and $X \backslash \backslash J$ are subcomplexes of $X$.

\medskip \noindent
\textbf{Theorem 1.2.2:} [Hag08, Corollary 2.16] \textit{Let $X$ be a CAT(0) cube complex and $J$ a hyperplane. Then $X \backslash \backslash J$ has exactly two connected components, which are combinatorially convex. The components of $\partial N(J)$ are also combinatorially convex.}

\medskip \noindent
The two connected components of $X \backslash \backslash J$ will be refered to as the \textit{halfspaces} associated to the hyperplane $J$.

In fact, the combinatorial metric and the hyperplanes are linked together: it can be proved that the combinatorial distance between two vertices corresponds exactly to the number of hyperplanes separating them [Hag08, Theorem 2.7].

Later, the following lemma will be useful. Roughly speaking, it proves that the natural \textit{combinatorial projection} is well-defined.

\medskip \noindent
\textbf{Lemme 1.2.3 :} \textit{Let $X$ be a CAT(0) cube complex, $C \subset X$ be a combinatorially convex subspace and $x \in X \backslash C$ be a vertex. Then there exists a unique vertex $y \in C$ minimizing the distance to $x$. Moreover, for any vertex of $C$, there exists a combinatorial geodesic from it to $x$ passing through $y$.}

\medskip \noindent
\textbf{Proof.} Because $d_c$ takes integer values, there clearly exists at least one vertex of $C$ minimizing the distance to $x$. Let $y_1,y_2$ be two vertices of $C$ minimizing the distance to $x$. For convenience, let $m$ denote the median point $m(x,y_1,y_2)$. Then, there exist combinatorial geodesics $[x,y_1]$, $[y_1,y_2]$, $[y_2,x]$ such that

\begin{center}
$\{ m \} = [x,y_1] \cap [y_1,y_2] \cap [y_2,x]$.
\end{center}

\noindent
In particular, $d(x,m) \leq \min (d(x,y_1),d(x,y_2))$. Now, $m \in [y_1,y_2] \subset C$ because $C$ is combinatorially convex, and since $y_1$ and $y_2$ minimize the distance to $x$, $d(x,y_1),d(x,y_2) \leq d(x,m)$. Consequently, $y_1=m=y_2$. Therefore, we have just proved that there exists at most one vertex of $C$ minimizing the distance to $x$.

\medskip \noindent
Now, let $z \in C$ be a vertex. Exactly the same argument as above gives $y=m(x,y,z)$, so there exists a combinatorial geodesic from $x$ to $z$ passing through $y$. $\square$

\medskip \noindent
\textbf{1.3. Squier and Farley complexes.} Let $\mathcal{P}= \langle \Sigma \mid \mathcal{R} \rangle$ be a semigroup presentation and $w \in \Sigma^+$ be a base word. One way of obtaining information about the diagram group $D(\mathcal{P},w)$ is to describe it as the fundamental group of a cube complex.

More precisely, we define the \textit{Squier complex} $S(\mathcal{P})$ as the cube complex whose vertices are the words in $\Sigma^+$; whose (oriented) edges can be written as $(a,u \to v,b)$, where $u=v$ or $v=u$ belongs to $\mathcal{R}$, linking the vertices $aub$ and $avb$; and whose $n$-cubes similarly can be written as $(a_1,u_1 \to v_1, \ldots, a_n, u_n \to v_n, a_{n+1})$, spanned by the set of vertices $\{ a_1w_1 \cdots a_nw_na_{n+1} \mid w_i=v_i \ \text{or} \ u_i \}$.

\medskip \noindent
\textbf{Theorem 1.3.1:} [GS97, Theorem 6.1] \textit{$D(\mathcal{P},w) \simeq \pi_1( S(\mathcal{P}),w)$.}

\medskip \noindent
For convenience, $S(\mathcal{P},w)$ will denote the connected component of $S(\mathcal{P})$ containing $w$. Notice that two words $w_1,w_2 \in \Sigma^+$ are equal modulo $\mathcal{P}$ if and only if they belong to the same connected component of $S(\mathcal{P})$. Therefore, a consequence of theorem 1.3.1 is:

\medskip \noindent
\textbf{Corollary 1.3.2:} \textit{If $w_1,w_2 \in \Sigma^+$ are equal modulo $\mathcal{P}$, then $D(\mathcal{P},w_1) \simeq D(\mathcal{P},w_2)$.}

\medskip \noindent
It can be proved that $S(\mathcal{P},w)$ is nonpositively curved so that its universal cover is CAT(0). In [Far99], Farley gives a construction of this cover.

A semigroup diagram is \textit{thin} whenever it can be written as a sum of atomic diagrams. We define the \textit{Farley complex} $X(\mathcal{P},w)$ as the cube complex whose vertices are the reduced semigroup diagrams $\Delta$ over $\mathcal{P}$ satisfying $\mathrm{top}(\Delta)=w$, and whose $n$-cubes are spanned by the vertices $\{ \Delta \cdot P \mid P \leq \Gamma \}$ for some vertex $\Delta$ and some thin diagram $\Gamma$ with $n$ cells. In particular, two diagrams $\Delta_1$ and $\Delta_2$ are linked by an edge if and only if there exists an atomic diagram $A$ such that $\Delta_1= \Delta_2 \cdot A$.

\medskip \noindent
\textbf{Theorem 1.3.3:} [Far99, Theorem 3.13] \textit{$X(\mathcal{P},w)$ is a CAT(0) cube complex. Moreover it is complete, i.e., there is no increasing sequence of cubes in $X(\mathcal{P},w)$.}

\medskip \noindent
There is a natural action of $D(\mathcal{P},w)$ on $X(\mathcal{P},w)$, namely $(g, \Delta) \mapsto g \cdot \Delta$. Then

\medskip \noindent
\textbf{Proposition 1.3.4:} [Far99, Theorem 3.13] \textit{The action $D(\mathcal{P},w) \curvearrowright X(\mathcal{P},w)$ is free. Moreover, it is properly discontinuous if $\mathcal{P}$ is a finite presentation, and it is cocompact if and only if the class $[w]_{\mathcal{P}}$ of words equal to $w$ modulo $\mathcal{P}$ is finite.}

\medskip \noindent
To conclude, we notice that the map $\Delta \mapsto \mathrm{bot}(\Delta)$ induces the universal covering $X(\mathcal{P},w) \to S(\mathcal{P},w)$ and that the action of $\pi_1(S(\mathcal{P},w))$ on $X(\mathcal{P},w)$ coincides with the natural action of $D(\mathcal{P},w)$. Precisely:

\medskip \noindent
\textbf{Lemma 1.3.5:} \textit{The map $\Delta \mapsto \mathrm{bot}(\Delta)$ induces a cellular isomorphism from the quotient $X(\mathcal{P},w)/D(\mathcal{P},w)$ to $S(\mathcal{P},w)$.}

\medskip \noindent
\textbf{Proof.} According to [Far99, Proposition 3.8], two diagrams $\Delta_1, \Delta_2 \in X(\mathcal{P},w)$ belong to the same orbit under the action of $D(\mathcal{P},w)$ if and only if $\mathrm{bot}(\Delta_1)= \mathrm{bot}(\Delta_2)$. Consequently, the map $\Delta \mapsto \mathrm{bot}(\Delta)$ induces a bijection $\varphi$ from the vertices of $X(\mathcal{P},w)/D(\mathcal{P},w)$ onto the vertices of $S(\mathcal{P},w)$.

\medskip \noindent
Now, from the descriptions of the cubes in Squier and Farley complexes, it is clear that $\varphi$ sends the vertices of an $n$-cube in $X(\mathcal{P},w)$ onto the vertices of an $n$-cube in $S(\mathcal{P},w)$; conversely, the vertices of an $n$-cube in $S(\mathcal{P},w)$ are the images of the vertices of an $n$-cube in $X(\mathcal{P},w)$. Therefore, $\varphi$ induces a cellular isomorphism from $X(\mathcal{P},w)/D(\mathcal{P},w)$ to $S(\mathcal{P},w)$. $\square$

\section{Hyperplanes of Farley complexes}

\noindent
From now on, let us fix a semigroup presentation $\mathcal{P}= \langle \Sigma \mid \mathcal{R} \rangle$ with a base word $w \in \Sigma^+$. Furthermore, for any hyperplane $J$ of $X( \mathcal{P},w)$, we introduce the following notation:
\begin{itemize}
	\item[$\bullet$] $J^+$ is the halfspace associated to $J$ not containing $\epsilon(w)$,
	\item[$\bullet$] $J^-$ is the halfspace associated to $J$ containing $\epsilon(w)$,
	\item[$\bullet$] $\partial_{\pm} J$ is the intersection $\partial N(J) \cap J^{\pm}$.
\end{itemize}
\noindent
Finally, for any diagram $A$, let $\#A$ denote its number of cells.

\medskip \noindent
\textbf{Definition 2.1:} A diagram $\Delta$ is \textit{minimal} if its maximal thin suffix $F$ has exactly one cell. (The existence of a maximal thin suffix is given by [Far99, Lemma 2.3].)

\medskip \noindent
In the following, $\overline{\Delta}$ will denote the diagram $\Delta \cdot F^{-1}$, obtained from $\Delta$ by removing the suffix $F$. The following result uses minimal diagrams to describe hyperplanes in Farley complexes.

\medskip \noindent
\textbf{Proposition 2.2:} \textit{Let $J$ be a hyperplane of $X(\mathcal{P},w)$. Then there exists a unique minimal diagram $\Delta$ such that $J^+= \{ D \ \text{diagram} \mid \Delta \leq D \}$. Conversely, if $\Delta$ is a minimal diagram and $J$ the hyperplane dual to the edge $[\Delta, \overline{\Delta}]$, then $J^+ = \{ D \ \text{diagram} \mid \Delta \leq D \}$.}

\medskip \noindent
To prove the proposition, two lemmas are needed:

\medskip \noindent
\textbf{Lemme 2.3:} \textit{Let $\Delta$ be a reduced diagram. If $\Delta = A_1 \circ \cdots \circ A_n$ where each $A_i$ is an atomic diagram, then the path}

\begin{center}
$\epsilon(w), \ A_1, \ A_1 \circ A_2, \ldots, \ A_1 \circ \cdots \circ A_n= \Delta$
\end{center}

\noindent
\textit{defines a combinatorial geodesic from $\epsilon(w)$ to $\Delta$. Furthermore, every combinatorial geodesic from $\epsilon(w)$ to $\Delta$ has this form.}

\medskip \noindent
\textbf{Proof.} In $X(\mathcal{P},w)$, if two vertices $B_1$ and $B_2$ are adjacent, then there exists an atomic diagram $D$ such that $B_1= B_2 \cdot D$. In particular, $|\# B_1 - \# B_2|=1$, so the combinatorial distance between $\epsilon(w)$ and $\Delta$ is bounded from below by $\# \Delta$, so that every path as mentionned in our lemma is a combinatorial geodesic, since its length is exactly the number of cells of $\Delta$.

\medskip \noindent
Conversely, let $(D_0,\ldots,D_n)$ be a combinatorial geodesic from $\epsilon(w)$ to $\Delta$. We know that, for every $0 \leq i \leq n-1$, there exists an atomic diagram $\Delta_i$ such that $D_{i+1} = D_i \cdot \Delta_i$. In particular,

\begin{center}
$\Delta = \Delta_{n-1} \cdot \Delta_{n-2} \cdots \Delta_1 \cdot \epsilon(w)= \Delta_{n-1} \cdot \Delta_{n-2} \cdots \Delta_1$.
\end{center}

\noindent
From what we have said above, the path $(D_0,\ldots, D_n)$ must have its length equal to $\# \Delta$. But the product above has at most $n$ cells, so it cannot contain dipoles hence

\begin{center}
$\Delta = \Delta_{n-1} \circ \Delta_{n-2} \circ \cdots \circ \Delta_1$.
\end{center}

\noindent
The proof is complete. $\square$

\medskip \noindent
\textbf{Corollary 2.4:} \textit{Let $A,B \in X(\mathcal{P},w)$ be two reduced diagrams. Then we have $d_c(A,B)= \# (A^{-1} \cdot B)$. In particular, $d_c(\epsilon(w),A)= \# (A)$.}

\medskip \noindent
\textbf{Proof.} The previous lemma clearly implies that the combinatorial distance from the trivial diagram corresponds to the number of cells. On the other hand, if
\begin{center}
$A, \ A \cdot \Delta_1, \ldots, \ A \cdot (\Delta_1 \circ \cdots \circ \Delta_n)=B$ 
\end{center}
is a combinatorial geodesic in $X(\mathcal{P},w)$ (where the $\Delta_i$ are atomic), then
\begin{center}
$\epsilon(\mathrm{bot}(A)), \ \Delta_1, \ldots, \ \Delta_1 \circ \cdots \circ \Delta_n=A^{-1} \cdot B$
\end{center}
defines a combinatorial geodesic in $X(\mathcal{P}, \mathrm{bot}(A))$. Thus
\begin{center}
$d_c(A,B)=d_c(\epsilon(w), A^{-1} \cdot B)= \# ( A^{-1} \cdot B)$. 
\end{center}
\noindent
The proof is complete. $\square$ 

\medskip \noindent
\textbf{Lemma 2.5:} \textit{Let $J$ be a hyperplane of $X(\mathcal{P},w)$, $\overline{\Delta}$ a vertex of $J^-$ minimizing the distance to $\epsilon(w)$ and $\Delta$ be the vertex adjacent to $\overline{\Delta}$ so that $J$ is dual to the edge $[\Delta,\overline{\Delta}]$. Then there exist $x,a,y \in \Sigma^+$ and a $(a,\ast)$-diagram $\Pi$ such that $\mathrm{bot}(\overline{\Delta})=xay$ and $\Delta= \overline{\Delta} \circ (\epsilon(x)+\Pi + \epsilon(y))$. Furthermore,}
\begin{center}
$\partial_- J = \left\{ \overline{\Delta} \circ (X+ \epsilon(a)+Y) \mid X \ \text{is an} \ (x,\ast)- \text{diagram}, \ Y \ \text{is a} \ (y,\ast)- \text{diagram} \right\}$
\end{center}
\noindent
\textit{and}
\begin{center}
$\partial_+ J = \left\{ \overline{\Delta} \circ (X+ \Pi+Y) \mid X \ \text{is an} \ (x,\ast)- \text{diagram}, \ Y \ \text{is a} \ (y,\ast)- \text{diagram} \right\}$
\end{center}

\noindent
\textbf{Proof.} Let $\Gamma$ be a diagram of the form $\overline{\Delta} \circ (X+ \epsilon(a)+Y)$, and write $X+ \epsilon(a)+Y$ as a concatenation $A_1 \circ \cdots \circ A_n$ of atomic diagrams. Because $\epsilon(x)+ \Pi + \epsilon(y)$ commutes with each $A_i$, the following configuration exists:
\begin{center}
\includegraphics[scale=0.5]{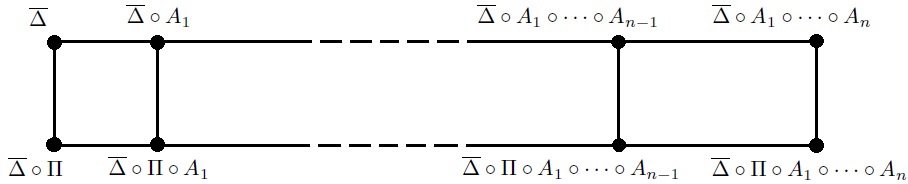}
\end{center}
\noindent
Consequently, $\overline{\Delta} \circ A_1 \circ \cdots \circ A_n \in \partial_-J$. We have thus proved the inclusions
\begin{center}
$\partial_- J \supset \{ \overline{\Delta} \circ (X+ \epsilon(a)+Y) \mid X \ (x,\ast)- \text{diagram}, \ Y \ (y,\ast)- \text{diagram} \}$
\end{center}
and
\begin{center}
$\partial_+ J \supset \{ \overline{\Delta} \circ (X+ \Pi+Y) \mid X \ (x,\ast)- \text{diagram}, \ Y \ (y,\ast)- \text{diagram} \}$
\end{center}
Conversely, let $\overline{\Gamma} \in \partial_- J$ be a vertex whose symmetric point with respect to the hyperplane $J$ will be denoted by $\Gamma$. We want to prove by induction on $d_c(\overline{\Delta},\overline{\Gamma})$ that $\overline{\Gamma}$ can be written as $\overline{\Delta} \circ (X+ \epsilon(a)+Y)$ and $\Gamma$ as $\overline{\Delta} \circ (X+ \Pi+Y)$, for some $(x,\ast)$-diagram $X$ and $(y,\ast)$-diagram $Y$. Of course, if $d_c(\overline{\Delta},\overline{\Gamma})=0$, there is nothing to prove since $\overline{\Gamma}= \overline{\Delta}$ and $\Gamma= \Delta$. So suppose that $d_c(\overline{\Delta}, \overline{\Gamma})=n+1$. According to Lemma 1.2.3, there exists a combinatorial geodesic from $\epsilon(w)$ to $\overline{\Gamma}$ passing through $\overline{\Delta}$. By convexity, the subpath from $\overline{\Delta}$ to $\overline{\Gamma}$ is included into $\partial_-J$; let $\overline{D}$ denote the next to last vertex. Thus, we have the configuration
\begin{center}
\includegraphics[scale=0.5]{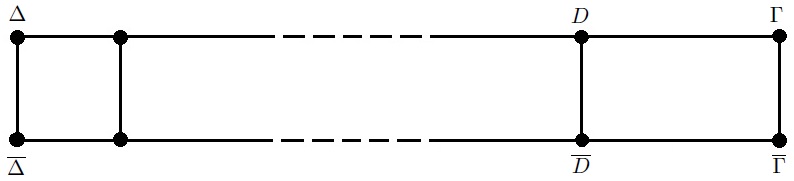}
\end{center}
\noindent
By the induction hypothesis, $\overline{D}= \overline{\Delta} \circ ( X+ \epsilon(a)+Y)$ and $D= \overline{\Delta} \circ ( X+ \Pi+Y)$ for some $(x,\ast)$-diagram $X$ and $(y,\ast)$-diagram $Y$. To the square induced by $D$, $\overline{D}$, $\Gamma$ and $\overline{\Gamma}$ is associated a thin diagram $F$, that we write as $F_1+F_2$ where $F_1$ and $F_2$ have exactly one cell. In particular, $F_1$ or $F_2$ has to contain $\Pi$ as a subdiagram, say $F_1$. Then, because $\overline{D}$ is the diagram with the minimal number of cells in the considered square, we may write $\mathrm{top}(F_1)=f_1$, $\mathrm{top}(F_2)=f_2$, and $\mathrm{bot}(Y)=gf_2h$, such that
\begin{center}
$\overline{\Gamma} = \overline{D} \circ ( \epsilon(f_1)+F_2) = \overline{D} \circ (X+ \epsilon(a)+Y) \circ (\epsilon(f_1)+F_2) = \overline{\Delta} \circ (X+ \epsilon(a)+Y')$,
\end{center}
\noindent
with $Y'= Y \circ (\epsilon(g)+F_2+ \epsilon(h))$; we also deduce that $\Gamma= \overline{\Delta} \circ (X+ \Pi+Y')$. Thus, we have just proved the inclusions
\begin{center}
$\partial_- J \subset \{ \overline{\Delta} \circ (X+ \epsilon(a)+Y) \mid X \ (x,\ast)- \text{diagram}, \ Y \ (y,\ast)- \text{diagram} \}$
\end{center}
\noindent
and
\begin{center}
$\partial_+ J \subset \{ \overline{\Delta} \circ (X+ \Pi+Y) \mid X \ (x,\ast)- \text{diagram}, \ Y \ (y,\ast)- \text{diagram} \}$.
\end{center}
\noindent
The proof is complete. $\square$

\medskip \noindent
\textbf{Proof of Proposition 2.2.} Let $J$ be a hyperplane. Let $\Delta$ denote the unique vertex of $J^+$ minimizing the distance to $\epsilon(w)$, whose existence is guaranted by Lemma 1.2.3. According to Lemma 1.2.3, for every vertex of $J^+$, there exists a combinatorial geodesic linking it with $\epsilon(w)$ and passing through $\Delta$, so we are able to deduce from Lemma 2.3 that $J^+ \subset \{ D \ \text{diagram} \mid \Delta \leq D \}$. Conversely, if $D$ is a diagram with $\Delta$ as a prefix, then Lemma 2.3 implies there exists a combinatorial geodesic from $\epsilon(w)$ to $D$ passing through $\Delta$; because $\Delta \in J^+$, we deduce that $D \in J^+$ by convexity. Therefore, we have $J^+ = \{ D \ \text{diagram} \mid \Delta \leq D \}$. 

\medskip \noindent
Now, we want to prove $\Delta$ is minimal. If it were not the case, we could write $\Delta = D \circ F$ with $F$ a thin diagram having exactly two cells. In this case, by construction of $X(\mathcal{P},w)$, $\Delta$ would be the vertex of a square $K$ with a minimal number of cells among the other three vertices, i.e., $\Delta$ would be farther from $\epsilon(w)$ than the other vertices (Corollary 2.4), that would imply that $J$ intersects $K$ along two adjacent edges in order to separate $\Delta$ from the other vertices of $K$, a contradiction because a hyperplane does not self-intersect in a CAT(0) cube complex [Sag95, Theorem 4.10].

\medskip \noindent
For the uniqueness of $\Delta$, it is sufficient to verify that $\Delta$ is the only minimal diagram in $\partial_+J$. According to Lemma 2.5, every diagram $\Gamma$ of $\partial_+J$ can be written as $\overline{\Delta} \circ (X+ \Pi + Y)$. In particular, the maximal thin suffix of $\Gamma$ contains $\Pi$, and a cell of $X$ and $Y$ if they are non-trivial. Consequently, if $\Gamma$ is minimal, necessarily $X$ and $Y$ must be trivial, hence 
\begin{center}
$\Gamma= \overline{\Delta} \circ (\epsilon(x)+ \Pi + \epsilon(y))= \Delta$. 
\end{center}
\noindent
Thus, the uniqueness is proved. $\square$

\medskip \noindent
\textbf{Proposition 2.6 :} \textit{Let $J$ be a hyperplane of $X(\mathcal{P},w)$ with $\Delta$ its associated minimal diagram. Let us write $\mathrm{bot}(\Delta)=aub$ so that the maximal thin suffix of $\Delta$ is of the form $\epsilon(a)+ \Pi + \epsilon(b)$, where $\Pi$ corresponds to a relation $u=v \in \mathcal{R}$. We will suppose that there is no word $\xi \in \Sigma^+$ such that $a=au \xi$ modulo $\mathcal{P}$. Then }
\begin{center}
$\mathrm{stab}(J)= \overline{\Delta} \cdot \left( D( \mathcal{P},a) \times D( \mathcal{P},b) \right) \cdot \overline{\Delta}^{-1}$,
\end{center}
\noindent
\textit{where $D( \mathcal{P},a) \times D( \mathcal{P},b)$ is the canonical subgroup of $D(\mathcal{P},aub)$ associated to the embedding $(A,B) \mapsto A+ \epsilon(u)+B$. }

\medskip \noindent
\textbf{Proof.} Let $\Gamma \in \mathrm{stab}(J)$. According to Lemma 3.1, independently proved in the next section, the action $D(\mathcal{P},w) \curvearrowright X(\mathcal{P},w)$ does not inverse hyperplanes, so $\Gamma$ stabilizes $\partial_-J$. According to Lemma 2.5, 
\begin{center}
$\partial_- J= \{ \overline{\Delta} \circ (A+ \epsilon(u)+B) \mid A \ (a,\ast)- \text{diagram}, \ B \ (b,\ast)-\text{diagram} \}$.
\end{center}
\noindent
It is worth noticing that $\Gamma$ stabilizes $\partial_-J$ if and only if $\overline{\Delta} \cdot \Gamma \cdot \overline{\Delta}^{-1}$ stabilizes
\begin{center}
$\{ A+ \epsilon(u)+B \mid A \ (a,\ast)- \text{diagram}, \ B \ (b,\ast)-\text{diagram} \}$.
\end{center}
\noindent
Thus, up to a conjugation by $\overline{\Delta}$, we may suppose without loss of generality that $\overline{\Delta}=\epsilon(w)$. In particular, $\epsilon(w) \in \partial_-J$, so $\Gamma = \Gamma \cdot \epsilon(w) \in \partial_-J$ because $\Gamma \in \mathrm{stab}(\partial_-J)$. Therefore, it is possible to write 
\begin{center}
$\Gamma = A+ \epsilon(u)+B$
\end{center}
\noindent
where $A$ is a $(a,\ast)$-diagram and $B$ a $(b,\ast)$-diagram. For convenience, let $\mathrm{bot}(A)= \alpha$ and $\mathrm{bot}(B)= \beta$. Because $\Gamma$ is a spherical diagram,
\begin{center}
$aub= \mathrm{top}(\Gamma) = \mathrm{bot}(\Gamma)= \alpha u \beta$. \hspace{2cm}         $(1)$
\end{center}
\noindent
Without loss of generality, we may suppose that $|\alpha | \geq |a |$ (where $| \cdot |$ is the length of a word written over $\Sigma$); otherwise, just replace $\Gamma$ with $\Gamma^{-1}$. If $|\alpha| > |a|$, we deduce from $(1)$ that $au$ is a prefix of $\alpha$, i.e., there exists a word $\xi \in \Sigma^+$ such that $\alpha = au \xi$ in $\Sigma^+$. So $A$ is a $(a,au\xi)$-diagram, hence $a=au \xi$ modulo $\mathcal{P}$, a relation forbidden by the hypotheses. So $|\alpha|=|a|$, i.e., $\alpha=a$ according to $(1)$. The relation $(1)$ then becomes $aub=au \beta$, hence $\beta=b$.

\medskip \noindent
Thus, $A$ and $B$ are respectively a $(a,a)$-diagram and a $(b,b)$-diagram, i.e., $\Gamma$ belongs to the canonical subgroup $D(\mathcal{P},a) \times D( \mathcal{P},b)$.

\medskip \noindent
Conversely, thanks to Lemma 2.5, it is clear that any element of $D(\mathcal{P},a) \times D( \mathcal{P},b)$ stabilizes $\partial_+J$ and $\partial_- J$, and in particular the hyperplane $J$ itself. $\square$

\section{Hyperplanes of Squier complexes}

\noindent
From now on, fix a semigroup presentation $\mathcal{P}= \langle \Sigma \mid \mathcal{R} \rangle$ with a base word $w \in \Sigma^+$.

\medskip \noindent
\textbf{Lemma 3.1:} \textit{Two (oriented) edges $(p, a \to b , q)$ and $(r, c \to d ,s)$ of $S(\mathcal{P})$ are dual to the same hyperplane if and only if $a=c$, $b=d$ in $\Sigma^+$ and $p=r$, $q=s$ modulo $\mathcal{P}$.}

\medskip \noindent
\textbf{Proof.} Suppose first that the edges $(p,a \to b,q)$ and $(r,c \to d,s)$ are dual to the same hyperplane. We show by induction on the length $\ell$ of a path of parallel edges that $(r,c \to d,s)$ has the desired form. If $\ell=0$, we have $(p,a \to b,q)=(r,c \to d,s)$ and there is nothing to prove. So suppose that $\ell=n+1$. By the induction hypothesis, the $n$-th edge of our path has the form $(p',a \to b,q')$, with $p'=p$, $q'=q$ modulo $\mathcal{P}$, and is parallel to $(r,c \to d,s)$. Therefore, these two edges belong to a same square 
\begin{center}
$(x,u \to v,y,a \to b,q')$ or $(p',a \to b, x,u \to v, y)$.
\end{center}
\noindent
Consequently,
\begin{center}
$(r, c \to d,s)= (xwy,a \to b,q')$ or $(p',a \to b,xwy)$
\end{center}
\noindent
with $w=u$ or $v$. Therefore, the edge $(r,c \to d,s)$ has the given form.

\medskip \noindent
Conversely, we want to show that the edges $(p,a \to b,q)$ and $(r, a \to b,s)$ are dual to the same hyperplane whenever $p=r$ and $q=s$ modulo $\mathcal{P}$. To do this, let us consider two sequences
\begin{center}
$p=x_1, \ x_2, \ldots, \ x_n=r$ and $q=y_1, \ y_2, \ldots, \ y_m=s$
\end{center}
\noindent
where $x_i=x_{i+1}$ and $y_i=y_{i+1}$ are relations of $\mathcal{R}$. Then, we have the configuration
\begin{center}
\includegraphics[scale=0.55]{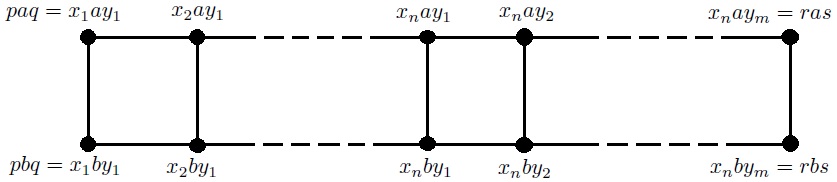}
\end{center}
\noindent
Therefore, our two edges are dual to the same hyperplane. $\square$

\medskip \noindent
\textbf{Lemma 3.2:} \textit{Each hyperplane of $S(\mathcal{P})$ is 2-sided.}

\medskip \noindent
\textbf{Proof.} It was not the case, there would exist two oriented edges $(a,x \to y,b)$ et $(a,y \to x,b)$ dual to the same hyperplane. According to Lemma 3.1, this is impossible. $\square$

\medskip \noindent
\textbf{Lemma 3.3:} \textit{A hyperplane of $S(\mathcal{P})$ self-intersects if and only if it is dual to an edge $(a,p \to q, bpc)$ where $a=apb$, $c=bpc$ modulo $\mathcal{P}$.}

\medskip \noindent
\textbf{Proof.} Let $J$ be a hyperplane dual to an edge $(a,p \to q,bpc)$ with $a=apb$ and $c=bpc$ modulo $\mathcal{P}$. Then, according to Lemma 3.1, the edges $(a,p \to q, bpc)$ and $(apb, p \to q,c)$ are dual to the same hyperplane; furthermore, they are adjacent edges in the square $(a,p \to q,b,p \to q,c)$. Consequently, $J$ self-intersects.

\medskip \noindent
Conversely, suppose that $J$ self-intersects. Because the relation associated to an edge dual to a fixed hyperplane is always the same (Lemma 3.1), $J$ self-intersects inside a square $(a,p \to q,b,p \to q,c)$. In particular, the edges $(a,p \to q,bpc)$ and $(apb,p \to q,c)$ are dual to our hyperplane $J$, hence the equalities $a=apb$ and $c=bpc$ modulo $\mathcal{P}$ according to Lemma 3.1. $\square$

\medskip \noindent
\textbf{Remark 3.4:} Lemmas 3.2 and 3.3 allow us to know whether the Squier complex $S(\mathcal{P},w)$ contains one-sided or self-intersecting hyperplanes directly from the semigroup presentation $\mathcal{P}$ and the base word $w$. Such hyperplanes belong to the four \textit{pathological configurations} identified by Haglund and Wise in [HW08]. Going further, it is possible to determine exactly when $S(\mathcal{P},w)$ is \textit{special} from $\mathcal{P}$ and $w$. For instance, we get the following simple criterion:

\medskip \noindent
\textit{Let $\mathcal{P}= \langle \Sigma \mid \mathcal{R} \rangle$ be a semigroup presentation and $w \in \Sigma^+$ a base word. If there are no words $a,b,p \in \Sigma^+$ satisfying the equalities $w=ab$, $a=ap$ and $b=pb$ modulo $\mathcal{P}$, then the Squier complex $S(\mathcal{P},w)$ is special. }

\medskip \noindent
In particular, it allows us to embed the associated diagram group into a (possibly non-finitely-generated) right-angled Artin group, giving linearity and residual finiteness of the group. More information will be available in a forthcoming paper.

\section{Main theorem}

\noindent
We prove below our main theorem.

\medskip \noindent
\textbf{Theorem 4.1:} \textit{Let $G$ be a diagram group. Then $G$ is free if and only if it does not contain any subgroup isomorphic to $\mathbb{Z}^2$.}

\medskip \noindent
Essentially, Theorem 4.1 will be a consequence of the following lemma:

\medskip \noindent
\textbf{Lemma 4.2 :} \textit{Let $\mathcal{P}= \langle \Sigma \mid \mathcal{R} \rangle$ be a semigroup presentation, $w \in \Sigma^+$ a base word and $H$ a non-trivial finitely generated subgroup of $D(\mathcal{P},w)$. We suppose that $D(\mathcal{P},w)$ does not contain any subgroup isomorphic to $\mathbb{Z}^2$. Then there exists a hyperplane $J$ of $X(\mathcal{P},w)$, inducing an arboreal structure $T$ on $X(\mathcal{P},w)$ which yields a non-trivial splitting of $H$ over the trivial subgroup.}

\medskip \noindent
\textbf{Proof.} First, if $\Delta$ is the maximal common prefix of all the spherical diagrams of $H$, then we may suppose without loss of generality that $\Delta$ is trivial; otherwise, just conjugate by $\Delta$ (and replace the base word $w$ with $\mathrm{bot}(\Delta)$). Thus, there does not exist any non-trivial prefix common to all spherical diagrams of $H$.

\medskip \noindent
Let $u$ be the maximal prefix of $w$ satisfying $D(\mathcal{P},u)= \{1 \}$. Notice that $u \neq w$, since otherwise $D(\mathcal{P},w)$ would be trivial, which is forbidden by the hypotheses. Thus, we can write $w=uav$ where $a \in \Sigma$ is the letter just following $u$ and where $v \in \Sigma^+$ is a possibly empty suffix. We also notice that $D(\mathcal{P},v)= \{1 \}$ from 
\begin{center}
$D(\mathcal{P},ua) \times D(\mathcal{P},v) \hookrightarrow D(\mathcal{P},w)$, 
\end{center}
\noindent
from $D(\mathcal{P},ua) \neq \{1 \}$ and from the fact that $D(\mathcal{P},w)$ does not contain any subgroup isomorphic to $\mathbb{Z}^2$, since diagram groups are torsion-free [GS97, Theorem 15.11]. Let $\Delta_0 \in H$ be a non-trivial spherical diagram. Because $\Delta_0 \notin D(\mathcal{P},u) \cup D(\mathcal{P},v)$, $\Delta_0$ contains a cell $\Pi$ just below the letter $a$, viewed as a subpath of $\mathrm{top}(\Delta_0)$. Let $\Delta_m$ denote the smallest prefix of $\Delta_0$ containing $\Pi$: we claim it is a minimal diagram. Indeed, otherwise $\Delta_m$ would have a maximal thin suffix with at least two cells; removing one of them, different from $\Pi$, we would obtain a smaller prefix containing $\Pi$, contradicting the minimality of $\Delta_m$. Let $J$ denote the hyperplane associated to $\Delta_m$ (Proposition 2.2).

\medskip \noindent
For convenience, say that $\mathrm{bot}(\overline{\Delta_m})=efagh$ where $fag=k$ is the relation associated to $\Pi$, so that $\mathrm{bot}(\Delta_m)=ekh$. That is to say that $\Delta_m$ has the form

\begin{center}
\includegraphics[scale=0.6]{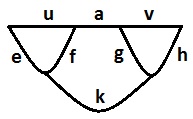}
\end{center}

\noindent
Here, it is possible that the cells associated to the relations $u \to ef$ and $v \to gh$ turn out to be trivial.

\medskip \noindent
We want to show that two translates of $J$ by $D(\mathcal{P},w)$ are either equal or disjoint. If it is not the case, the image $\overline{J}$ of $J$ in the Squier complex $S(\mathcal{P})$ self-intersects. According to Lemma 3.3, $\overline{J}$ has to be dual to an edge $(r,p \to q,spt)$ with $r=rps$ and $s=spt$ modulo $\mathcal{P}$. But we already know that $\overline{J}$ is dual to the edge $(e,fag \to k,h)$, hence, according to Lemma 3.1, $e=r$, $h=t$ modulo $\mathcal{P}$ and $p=fag$, $q=k$ in $\Sigma^+$. Consequently, using the embeddings associated to the canonical subgroups and Corollary 1.3.2,
\begin{center}
$D(\mathcal{P},ua) \hookrightarrow D(\mathcal{P},uag) \simeq D(P,efag)  \simeq D(P,rp)$,
\end{center}
\noindent
and then
\begin{center}
$D(\mathcal{P},rp) \hookrightarrow D(\mathcal{P},rps) \simeq D(\mathcal{P},r) \simeq D(\mathcal{P},e) \hookrightarrow D(\mathcal{P},ef) \simeq D(\mathcal{P},u)$,
\end{center}
\noindent
but it is impossible because $D(\mathcal{P},ua) \neq \{1 \}$ and $D(\mathcal{P},u)= \{1 \}$.

\medskip \noindent
Consequently, the orbit $H \cdot J$ induces an arboreal structure on $X(\mathcal{P},w)$, i.e., if we define the graph $T$ whose vertices are the components of $X(\mathcal{P},w) \backslash D(\mathcal{P},w) \cdot J$ and where two components are linked by an edge precisely when they are separated by exactly one translate of $J$, then $T$ is a tree; moreover $H$ acts naturally on $T$. According to Proposition 2.2, $J$ separates $\epsilon(w)$ and $\Delta_0$, so that $T$ has at least two vertices. 

\medskip \noindent
Notice that the stabilizer of $J$ is trivial. To see this, we will apply Proposition 2.6, i.e., we claim that there is no word $\xi \in \Sigma^+$ such that $e=efag \xi$ modulo $\mathcal{P}$. If it was the case, we would have
\begin{center}
$D(\mathcal{P},ua) \hookrightarrow D(\mathcal{P},efag \xi) \simeq D(\mathcal{P},e) \hookrightarrow D(\mathcal{P},u)$,
\end{center}
\noindent
which is impossible because $D(\mathcal{P},ua) \neq \{1 \}$ and $D(\mathcal{P},u)=\{1 \}$. Consequently, Proposition 2.6 gives
\begin{center}
$\mathrm{stab}(J) = \overline{\Delta_m} \cdot D(\mathcal{P},e) \times D ( \mathcal{P},h) \cdot \overline{\Delta_m}^{-1}= \{1 \}$,
\end{center}
\noindent
because $D(\mathcal{P},e) \hookrightarrow D(\mathcal{P},u)= \{ 1 \}$ and $D(\mathcal{P},h) \hookrightarrow D(\mathcal{P},v)= \{1 \}$.

\medskip \noindent
To sum up, we have a tree $T$ with at least two vertices, on which acts $H$ with only one orbit of edges whose stabilizer is trivial. Now, two cases may happen : either the image $\overline{J}$ of $J$ in the Squier complex separates $S(\mathcal{P},w)$ into exactly two connected components, or $S(\mathcal{P},w) \backslash \overline{J}$ is connected. In the second case, we deduce that $T$ has exactly one orbit of vertices, i.e., the splitting of $H$ associated to $T$ is an HNN extension over the trivial subgroup; of course, such a splitting cannot be trivial. 

\medskip \noindent
From now on, suppose that $\overline{J}$ separates $S(\mathcal{P},w)$ into two connected components. This allows us to classify the connected components of $X(\mathcal{P},w) \backslash H \cdot J$ into two families: those which are sent onto the component of $S(\mathcal{P},w) \backslash \overline{J}$ containing $w$, and the others; according to Lemma 1.3.5, it amounts to distinguish components containing or not spherical diagrams. Thus, $T$ has naturally a bipartite structure, moreover $H$-invariant because the set of spherical diagrams is itself $H$-invariant; in particular, we deduce that the action on $T$ is without edge-inversions, which we already knew from Lemma 3.2. Now, we have to verify that the action has no global fixed point. If such a fixed point existed, one of the two families defined above would be reduced to only one component of $X(\mathcal{P},w) \backslash H \cdot J$. Because $J$ separates $\epsilon(w)$ and $\Delta_0$, there exist at least two components containing spherical diagrams. From now on, let us suppose there exists only one component not containing spherical diagrams; from the bipartite structure of $T$, we deduce that $T$ is a \textit{star}, i.e., a vertex $s$ fixed by $H$ with adjacent vertices $s_i$ permuted by $H$. Therefore, $X(\mathcal{P},w)$ can be described as

\begin{center}
\includegraphics[scale=0.6]{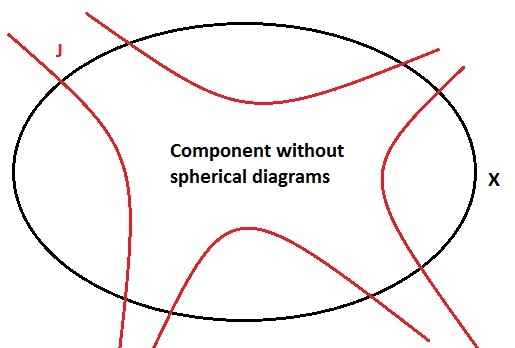}
\end{center}

\noindent
From the action $H \curvearrowright T$, we are able to write $H$ as the free product between the stabilizer of $s$, i.e., $H$ itself, and the stabilizer of one of the $s_i$. Therefore, since $H$ is finitely-generated, we deduce from Grushko's theorem that the stabilizer of a connected component of $X(\mathcal{P},w) \backslash H \cdot J$ containing at least one spherical diagram must be trivial. Now, we want to argue that this is impossible. Indeed, because the spherical diagrams of $H$ have no common non-trivial prefix, we deduce that there exists a non-trivial spherical diagram $\Delta_1 \in H$ such that $\Delta_m$ is not one of its prefix; according to Proposition 2.2, $J$ separates $\Delta_0$ and $\{ \epsilon(w), \Delta_1 \}$. Consequently, $\Delta_1 \in H \backslash \{1 \}$ belongs to the stabilizer of the connected component of $X(\mathcal{P},w) \backslash H \cdot J$ containing $\epsilon(w)$, so this stabilizer cannot be trivial. 

\medskip \noindent
We have just proved that the action $H \curvearrowright T$ has no global fixed point. Therefore, the associated splitting over the trivial subgroup is not trivial. $\square$

\medskip \noindent
\textbf{Proof of Theorem 4.1 :} Suppose by contradiction that $G$ is a diagram group not containing a subgroup isomorphic to $\mathbb{Z}^2$ but that is not free. Therefore, a non-trivial relation holds in $G$, and such a relation can be read in a non-trivial finitely-generated subgroup $H$. Furthermore, without loss of generality, $H$ may be supposed to be freely irreducible; otherwise, write $H$ as a free product of finitely-many freely irreducible groups $H_1 \ast \cdots \ast H_n$ and notice that at least one $H_i$ is not free. However, according to Lemma 4.2, $H$ must be either trivial or infinite cyclic, a contradiction. $\square$

\medskip \noindent
\textbf{Corollary 4.3 :} \textit{Let $G$ be a finitely-generated diagram group. Then $G$ is hyperbolic if and only if $G$ does not contain any subgroup isomorphic to $\mathbb{Z}^2$ if and only if $G$ is free.}

\section{Some examples}

\noindent
In this last section, we want to show how Theorem 4.1 can be used in practice to determine whether or not a given diagram group is free. First, we need a criterion to determine when a diagram group contains $\mathbb{Z}^2$. Essentially, it will be a consequence of the following result due to Guba and Sapir [GS99, Theorem 17]:

\medskip \noindent
\textbf{Proposition 5.1:} \textit{Let $\mathcal{P}= \langle \Sigma \mid \mathcal{R} \rangle$ be a semigroup presentation and let $G= D(\mathcal{P},w)$ for some $w \in \Sigma^+$. Suppose that $A_1, \ldots, A_m$ are spherical diagrams with base $w$ that pairwise commute in $G$. Then there exist a word $v=v_1 \cdots v_n$, spherical $(v_j,v_j)$-diagrams $\Delta_j$ for $1 \leq j \leq n$, integers $d_{ij}$ for $1 \leq i \leq m$ and $1 \leq j \leq n$, and some $(w,v)$-diagram $\Gamma$ such that
\begin{center}
$A_i = \Gamma \cdot \left( \Delta_1^{d_{i1}}+ \cdots + \Delta_n^{d_{in}} \right) \cdot \Gamma^{-1}$
\end{center}
\noindent
for all $1 \leq i \leq m$.} 

\medskip \noindent
As a corollary, we obtain a method to compute the \textit{algebraic dimension} of a diagram group, i.e., the higher possible rank of a free Abelian subgroup.

\medskip \noindent
\textbf{Corollary 5.2:} \textit{Let $\mathcal{P}= \langle \Sigma \mid \mathcal{R} \rangle$ be a semigroup presentation and $w \in \Sigma^+$ a base word. Then $D(\mathcal{P},w)$ contains a free Abelian subgroup of rank $n$ if and only if there exist words $w_1,\ldots,w_n \in \Sigma^+$ such that $w=w_1 \cdots w_n$ modulo $\mathcal{P}$ and $D(\mathcal{P},w_i) \neq \{1 \}$ for all $1 \leq i \leq n$.}

\medskip \noindent
\textbf{Proof.} If there exist words $w_1,\ldots,w_n \in \Sigma^+$ such that $w=w_1 \cdots w_n$ modulo $\mathcal{P}$ and $D(\mathcal{P},w_i) \neq \{1 \}$ for all $1 \leq i \leq n$, then
\begin{center}
$\mathbb{Z}^n = \mathbb{Z} \times \cdots \times \mathbb{Z} \hookrightarrow D(\mathcal{P},w_1) \times \cdots \times D(\mathcal{P},w_n) \hookrightarrow D(\mathcal{P},w)$.
\end{center}
\noindent
Conversely, suppose that $D(\mathcal{P},w)$ contains a free Abelian subgroup of rank $n$. Let $A_1, \ldots, A_n$ be a free basis of such a subgroup. According to Proposition 5.1, there exist a word $v=v_1 \cdots v_m$, spherical $(v_j,v_j)$-diagrams $\Delta_j$ for $1 \leq j \leq m$, integers $d_{ij}$ for $1 \leq i \leq n$ and $1 \leq j \leq m$, and some $(w,v)$-diagram $\Gamma$ such that
\begin{center}
$A_i= \Gamma \cdot \left( \Delta_1^{d_{i1}}+ \cdots + \Delta_n^{d_{in}} \right) \cdot \Gamma^{-1}$
\end{center}
\noindent
for all $1 \leq i \leq n$. For convenience, let $\Gamma_j$ denote the spherical diagram
\begin{center}
$\Gamma \cdot \left( \epsilon(v_1 \cdots v_{j-1}) + \Delta_j + \epsilon(v_{j+1} \cdots v_m) \right) \cdot \Gamma^{-1}$.
\end{center}
\noindent
Clearly, $\langle \Gamma_1, \ldots, \Gamma_m \rangle$ is a free Abelian subgroup of rank $k$ containing $\langle A_1,\ldots, A_n \rangle$, where $k$ is the number non-trivial $\Delta_i$. In particular, we deduce that $k \geq n$. For convenience, we will suppose that $\Delta_1, \ldots, \Delta_{k-1}$ are non-trivial in order to shorten the notation. The general case is proved exactly with the same argument. 

\medskip \noindent
Now, modulo $\mathcal{P}$, $w=w_1 \cdots w_n$ where $w_i=v_i$ for $1 \leq i \leq n-1$ and $w_n=v_n \cdots v_m$. Moreover, $D(\mathcal{P},w_i) \neq \{ 1 \}$ because $\Delta_i$ is a non-trivial spherical $(w_i,w_i)$-diagram for $1 \leq i \leq n-1$, and $\Delta_n+ \cdots + \Delta_m$ is a non-trivial spherical $(w_n,w_n)$-diagram. $\square$

\medskip \noindent
Therefore, combining Corollary 5.2 with Theorem 4.1, we deduce a method, often quite simple, to determine whether or not a given diagram group is free. Let us consider a few examples.

\medskip \noindent
\textbf{Example 5.3:} Let $\mathcal{P}= \langle a,b,c \mid ab=ba, ca=ac, bc=cb \rangle$ and $U(p,q,r)=D(\mathcal{P},a^pb^qc^r)$. In [GS97, Example 10.2], Guba and Sapir prove that $U(p,q,r)$ is free if $p \leq 1$, $q \leq 1$ or $r \leq 1$. In fact, it is sufficient to prove that $U(1,q,r)$ is free.

\medskip \noindent
If we want to break a word $w$ equal to $ab^qc^r$ modulo $\mathcal{P}$ into two pieces $w=w_1w_2$, then $w_1$ or $w_2$, say $w_1$, contains only the letters $b$ and $c$. Therefore, $D(\mathcal{P},w_1)$ is naturally a subgroup of $U(0,q,r)$. However, noticing that the associated Squier complex is simply connected, we deduce that $U(0,q,r)= \{1 \}$, hence $D(\mathcal{P},w_1)= \{ 1 \}$. Thus, the algebraic dimension of $U(1,q,r)$ is $1$; in particular, it is free.

\medskip \noindent
Using the presentation given by [GS97, Theorem 9.8], it is not difficult to find a free basis of $U(1,q,r)$; in particular, we notice that its rank is $qr$.

\medskip \noindent
\textbf{Example 5.4:} Let $\mathcal{P} = \langle a,b \mid ab=a, ab^2=a \rangle$. We claim that $D(\mathcal{P},a)$ is free.

\medskip \noindent
Indeed, the set of words equal to $a$ modulo $\mathcal{P}$ is exactly $\{ ab^n \mid n\geq 0 \}$, and if we want to break $ab^n$ into two pieces, we obtain something like $ab^n = (ab^m)b^{n-m}$ with $D(\mathcal{P},ab^m) \simeq D(\mathcal{P},a) \neq \{1 \}$ and $D(\mathcal{P},b^{n-m})= \{1 \}$ (since any relation of $\mathcal{P}$ contains the letter $a$). Therefore, the algebraic dimension of $D(\mathcal{P},a)$ is $1$; in particular, it is free. 

\medskip \noindent
Using the presentation given by [GS97, Theorem 9.8], it is not difficult to find a free basis of $D(\mathcal{P},a)$; in particular, we notice that it has countable rank.

\medskip \noindent
\textbf{Example 5.5:} Let $\mathcal{P}= \langle a,b,c \mid ap=a, pc=c, bp=b,pb=b \rangle$. We claim that the group $D(\mathcal{P},abc)$ is free.

\medskip \noindent
Indeed, the set of words equal to $abc$ modulo $\mathcal{P}$ is exactly $\{ap^nbp^mc \mid n,m \geq 0 \}$, and if we want to break $ap^nbp^mc$ into two pieces, we get either $(ap^k)(p^{n-k}bp^mc)$ or $(ap^nbp^k)(p^{m-k}c)$; because $D(\mathcal{P},ap^k) \simeq D(\mathcal{P},a) = \{1 \}$ (since $S(\mathcal{P},a)$ is simply connected), $D(\mathcal{P},p^{n-k}bp^mc) \simeq D(\mathcal{P},bc) \neq \{1 \}$, $D(\mathcal{P},ap^nbp^k) \simeq D(\mathcal{P},ab) \neq \{1 \}$ and $D(\mathcal{P},p^{m-k}c) \simeq D(\mathcal{P},c) = \{1 \}$ (since $S(\mathcal{P},c)$ is simply connected), we deduce that the algebraic dimension of $D(\mathcal{P},abc)$ is $1$, and that it is a free group.

\medskip \noindent
Using the presentation given by [GS97, Theorem 9.8], it is not difficult to find a free basis of $D(\mathcal{P},abc)$; in particular, we find that it is a free group of rank two.

\medskip \noindent
\textbf{Examples 5.6:} In fact, because many semigroup presentations give free diagram groups, the method above turns out to be really useful. For instance, the same argument can be applied to prove that the groups $D(\mathcal{P}_i,w_i)$ are free, where $\mathcal{P}_i$ and $w_i$ are given by:
\begin{center}
$\mathcal{P}_1 = \langle a,c,p,x \mid ap=a, pc=c, ax=a, xp=p \rangle$ and $w_1=ac$,
\end{center}
\begin{center}
$\mathcal{P}_2 = \langle a,b,p,q \mid ap=a, pb=b, q=p, aqqb=appb \rangle$ and $w_2=ab$,
\end{center}
\begin{center}
$\mathcal{P}_3 = \langle a,x,y \mid ax=a, ay=a, y=x \rangle$ and $w_3= a$.
\end{center}
The details are left to the reader as an exercise.

\section*{References}
\addcontentsline{toc}{section}{References}

\begin{itemize}
	\item[$\text{[CH09]}$] P.-E. Caprace and F. Haglund, \textit{On geometric flats in the CAT(0) realization of Coxeter groups and Tits buildings}, Canad. J. Math. Vol. 61 (4), 2009 pp. 740-761.
	\item[$\text{[CSS07]}$] J. Crisp, M. Sageev and M. Sapir, \text{Surface subgroups of right-angled Artin groups}, arXiv:0707.1144 (2007).
	\item[$\text{[Che00]}$] V. Chepoi, \textit{Graphs of some CAT(0) complexes}, Adv. in Appl. Math., 24(2): 125-179, 2000.
	\item[$\text{[Far00]}$] D. Farley, \textit{Finiteness and CAT(0) properties of diagram groups}, PhD thesis (2000).
	\item[$\text{[GS97]}$] V. Guba and M. Sapir, \textit{Diagram groups}, Mem. Amer. Math. Soc 130 (620) (1997) 1-117.
	\item[$\text{[GS99]}$] V. Guba and M. Sapir, \textit{On subgroups of R. Thompson's group $F$ and other diagram groups}, Mat. Sb. 190 (8) (1999) 3-16 (in Russian). English translation in: V.S. Guba, M.V. Sapir, Sbornik: Mathematics 190 (8) (1999) 1077-1130.
	\item[$\text{[GS06a]}$] V. Guba and M. Sapir, \textit{Diagram groups and directed 2-complexes: homotopy and homology}, Journal of Pure and Applied Algebra 205 (2006) 1-47.
	\item[$\text{[GS06b]}$] V. Guba and M. Sapir, \text{Diagram groups are totally orderable}, Journal of Pure and Applied Algebra 205 (2006) 48-73.
	\item[$\text{[Hag08]}$] F. Haglund, \textit{Finite-index subgroups of graph products}, Geom. Dedicata (2008) 135:167-209.
	\item[$\text{[HW08]}$] F. Haglund and D. Wise, \textit{Special cube complexes}, GAFA 17(5), 1551-1620 (2008).
	\item[$\text{[Kil94]}$] V. Kilibarda, \textit{On the algebra of semigroup diagrams}, PhD thesis (1994). 
	\item[$\text{[Sa95]}$] M. Sageev, \textit{Ends of Group Pairs and Non-Positively Curved Cube Complexes}, Proc. London Math. Soc. (3) 71(3), 585-617 (1995).
	\item[$\text{[SW11]}$] M. Sageev and D. Wise, \textit{Periodic flats in CAT(0) cube complexes}, Algebraic and Geometric Topology 11 (2011) 1793-1820.
	\item[$\text{[S07]}$] M. Sapir, \textit{Some group theory problems}, International Journal of Algebra and Computation, Vol. 15, Nos. 5 \& 6 (2007) 1189-1214. 
	\item[$\text{[Se80]}$] J.-P. Serre, \textit{Trees}, Springer, Berlin (1980).
\end{itemize}

\end{document}